\documentclass[10pt,a4paper]{article}

\usepackage[margin=2.5cm]{geometry}
\usepackage{amsmath,amsthm,amsfonts, amssymb, array, color, graphicx, float, subcaption}
\usepackage{pgfplots, pgfplotstable}
\pgfplotsset{compat=1.15}
\usepackage{tikz}
\usepackage[linesnumbered,ruled,vlined]{algorithm2e}
\usepackage[braket, qm]{qcircuit}

\theoremstyle{definition}

\numberwithin{equation}{section}
\numberwithin{figure}{section}
\usepackage{listings}
\usepackage{color} 
\definecolor{mygreen}{RGB}{28,172,0} 
\definecolor{mylilas}{RGB}{170,55,241}
\font\myfont=cmr12 at 20pt
\title{{\myfont A Bisection Method Like Algorithm for Approximating Extrema of a Continuous Function}}

\author{
	Fatih İdiz\footnote{ \texttt{Email:fatihidiz@iyte.edu.tr}}
}

\begin{document}
	
	\maketitle
	
\begin{abstract}
For a continuous function $f$ defined on a closed and bounded domain, there is at least one maximum and one minimum. First, we introduce some preliminaries which are necessary through the paper. We then present an algorithm, which is similar to the bisection method, to approximate those maximum and minimum values. We analyze the order of the convergence of the method and the error at the $k$-th step. Then we discuss the pros and cons of the method. Finally, we apply our method for some special classes of functions to obtain nicer results. At the end, we write a Matlab script which implements our algorithm.
\end{abstract}

	
	\section{Introduction}
For nonlinear equations, there is no general algorithm or method that computes an exact root of a nonlinear equation, unfortunately. However, in science and engineering, one usually encounters a situation in which she needs to compute a root of a nonlinear equation. For this reason, there are several methods that approximates a root of a nonlinear equation, numerically. One of those methods is the bisection method, which is alternatively called binary chopping, interval halving, or bracketing method.

Before introducing and proving the bisection method, we need to recall some concepts and theorems from calculus. \\ 
{\bf 1.1. Monotone Convergence Theorem(Ross[1]).} All bounded monotone sequences converge.
\begin{proof}
See [1, Theorem 10.2.]
\end{proof} 

\par \noindent
{\bf 1.2. The Squeeze(Sandwich) Theorem([2]).} Suppose that $g(x)\leq f(x) \leq h(x)$ for all $x$ in some open interval containing $c$, except possibly at $x=c$ itself. Suppose also that $\displaystyle{\lim_{x \to c} g(x)=\lim_{x \to c} h(x) = L}$. Then $\displaystyle{\lim_{x \to c} f(x) = L}$.
\begin{proof}
See [2, Theorem 4] 
\end{proof} 
\par \noindent
{\bf 1.3. Intermediate-Value Theorem(Lang[3]).} Let $f$ be a continuous function on a closed interval $[a,b]$. Let $\alpha = f(a)$ and $\beta = f(b)$. Let $\gamma$ be a number such that $\alpha < \gamma < \beta$ or $\beta < \gamma < \alpha$. Then there exists a number $c$, $a<c<b$, such that $f(c)=\gamma$.
\begin{proof}
See [3, Theorem 4.4.]
\end{proof} 
\par \noindent
{\bf 1.4. Extreme Value Theorem(Lang[4]).} Let $f$ be a continuous function on a closed interval $[a,b]$. Then there exists an element $c\in [a,b]$, such that $f(c)\geq f(x)$ for all $x\in [a,b]$. That is, $c$ is a maximum for $f$ on $[a,b]$. Also, there exists an element $d\in[a,b]$ such that $f(x)\geq f(d)$ for all $x\in[a,b]$. That is, $d$ is a minimum for $f$ on $[a,b]$.
\begin{proof}
See [4, Theorem 4.3.]
\end{proof} 
\par \noindent
{\bf 1.5. The Bisection Method(Kincaid, Cheney[5]).} If $f$ is a continuous function on a closed interval $[a,b]$ and if $f(a)f(b)<0$, then $f$ must have a zero in $(a,b)$. Since $f(a)f(b)<0$, the function $f$ changes sign on the interval $[a,b]$ and, therefore, it has at least one zero in the interval. This is a consequence of the Intermediate-Value Theorem. The Bisection Method exploits this idea in the following way. If $f(a)f(b)<0$, then we compute $\displaystyle{c=\frac{a+b}{2}}$ and test whether $f(a)f(c)<0$. If this is true, then $f$ has a zero in $[a,c]$. So we rename $c$ as $b$ and start again with the new interval $[a,b]$, which is half as large as the original interval. If $f(a)f(c)>0$,then $f(c)f(b)<0$, and in this case we rename $c$ as $a$. In either case, a new interval containing a zero of $f$ has been produced, and the process can be repeated. The Bisection Method finds one zero but not all the zeros in the interval $[a,b]$. Of course, if $f(a)f(c)=0$, then $f(c)=0$ and a zero has been found. However, it is quite unlikely that $f(c)$ is exactly $0$ in the computer because of roundoff errors. Thus, the stopping criterion should not be whether $f(c)=0$. A reasonable tolerance must be allowed depending on the situation. \\ 

{\bf An example Matlab code for solving $e^x-4x=0$ in the interval $[0,1]$ by using the bisection method:}
\begin{lstlisting}[frame=single]
% set f to be a mathematical function
f = @(x) exp(x)-4*x;
% plot the function
fplot(f, [0,1])
% now we aim to solve exp(x)-4x == 0 to 6 decimal places
epsilon = 1e-7;
a = 0.0; b = 1;
while abs(b - a)/2 > epsilon
  c = (a+b)/2;
  if f(b)*f(c) > 0
       b = c;
  else
       a = c;
  end
end
estRoot = (a+b)/2
% check that f(estRoot) is indeed small
display(f(estRoot))
\end{lstlisting}
\par \noindent
{\bf Theorem 1.6.} If $[a_0,b_0 ],[a_1,b_1 ],\dots,[a_n,b_n ],\dots$ denote the intervals in the bisection method, then the limits $\displaystyle{\lim_{n \to \infty} a_n}$ and $\displaystyle{\lim_{n \to \infty} b_n}$ exist, are equal, and represent a zero of  $f$. If $r=\displaystyle{\lim_{n \to \infty} c_n}$ and $c_n=\displaystyle{\frac{a_n+b_n}{2}}$, then
\begin{equation*}
    |r-c_n|\leq 2^{-(n+1)}(b_0-a_0)
\end{equation*}
\begin{proof}
See [5, pp. 74-78]
\end{proof}

\section{Construction of Our Algorithm}

Let $f$ be a continuous function on a closed interval $[a,b]$ and $\alpha=f(a)$ and $\beta=f(b)$. By the Extreme Value Theorem, $f$ has at least one maximum value and one minimum value. Therefore, $f$ is bounded on the interval $[a,b]$. Let $U$ be an upper bound and $L$ be a lower bound for $f$ on the interval $[a,b]$.That is, $L\leq f(x)\leq U$ for all $x\in[a,b]$. Let $f$ attain its maximum value at $c\in[a,b]$ and $M=f(c)$. Since $f$ is continuous on $[a,b]$, $f$ is also continuous on $[a,c]$ and $[c,b]$. Let $m$ be the minimum of $\alpha=f(a)$ and $\beta=f(b)$, that is, $m=\text{min}\{\alpha,\beta\}$. Then there are two cases: 1) $m=\alpha$ or 2) $m=\beta$. 
\par \noindent
{\bf Case 1.} $m=\alpha$ 

If $m=\alpha$, then for all $y\in(m,M)$, there is a number $x\in(a,c)$ such that $y=f(x)$. Thus, there is no number $y$ between $m$ and $M$ such that the equation $y=f(x)$ has no solution in the interval $(a,c)$. Since $f(x)\leq U$ for all $x\in[a,b]$, we can conclude that $\alpha=m\leq M=f(c)\leq U$. Now, let $N$ be any number between $\alpha=m$ and $M=f(c)$. Then $N\leq M=f(c)\leq U$. Let $\displaystyle{A=\frac{N+U}{2}}$. If the equation $f(x)=A$ has a solution in the interval $[a,b]$, then $A\leq M=f(c)\leq U$ and we rename $A$ as $N$. So we start again with the new interval $[N,U]$, which is half as large as the original interval. If the equation $f(x)=A$ has no solution in the interval $[a,b]$, then $A$ is a new upper bound for the function $f$ on the interval $[a,b]$ because we have shown that for all $y\in(m,M)$, there is a number $x\in(a,c)$ such that $y=f(x)$. Thus, if the equation $f(x)=A$ has no solution in the interval $[a,b]$, then $A$ must be a new upper bound which is less than $U$. Then we rename $A$ as $U$. In either case, a new interval containing a maximum value of $f$ has been produced, and the process can be repeated. By repeating this process, we obtain a smaller and smaller interval that contains the maximum value of $f$.
\par \noindent
{\bf Case 2.} $m=\beta$ 

When $m=\beta$, a similar discussion like in Case 1 $(m=\alpha)$ can be carried out. Thus, in either case, our algorithm converges to the maximum value of the function.
\begin{algorithm}[h]
	\SetKwInOut{Input}{Input}\SetKwInOut{Output}{Output}
	\Input{1) Function $f$ \\ 2) Tolerance(=Error) $\epsilon_{tol}$ \\ 3) Points $m$ and $u$ such that $m=f(x)$ has a solution $x\in [a,b]$ and $u=f(x)$ \\ has no solution in $[a,b]$ \\ 4) Maximum iterations $N_{MAX}$ to prevent infinite loop}
	\Output{Value which differs from the maximum of $f(x)$ in the interval $[a,b]$ by less than $\epsilon_{tol}$}
	\caption{Our Constructed Algorithm\label{our algorithm}}
	\BlankLine
	$N\leftarrow 1$ \\
	\While{$N \leq N_{MAX}$}{
				$c\leftarrow (u+m)/2$\;
				\If{$(u-m)/2 \leq \epsilon_{tol}$}{
				Output($c$)\; \par
				\textbf{Stop}
			}
			$N=N+1$\;
			\If{$c=f(x)$ has a solution in $[a,b]$}{
				$m\leftarrow c$\;
			}
			\Else{
				$u\leftarrow c$\;
			}
			}
			Output(``Maximum iterations reached without the desired tolerance. Input a bigger $N_{MAX}$")
\end{algorithm}

\section{Error Analysis}

Let $m_0$ be a number such that $m\leq m_0\leq M=f(c)$ and $u_0$ be an upper bound for the function $f$ on the interval $[a,b]$. Let $[m_n,u_n ]$ be the interval at the $n$-th step. Then
\begin{equation*} \label{eq1}
\begin{split}
m_0 & \leq m_1 \leq m_2 \leq \dots \leq u_0 \\
u_0 & \geq u_1 \geq u_2 \geq \dots \geq m_0
\end{split}
\end{equation*}
\begin{flalign}\label{eq:1}
& & & \ \ \ u_{n+1}-m_{n+1} = \frac{u_n-m_n}{2}  & & (n\geq 0)
\end{flalign}
Since the sequence $(m_n)$ is nondecreasing (i.e. monotone) and bounded, it converges by The Monotone Convergence Theorem. Likewise, $(u_n)$ converges as it is nonincreasing (i.e. monotone) and bounded. If we apply Equation [\ref{eq:1}], repeatedly, we find that
\begin{equation*}
    u_n-m_n=2^{-n}(u_0-m_0)
\end{equation*}
Thus,
\begin{equation*}
    \lim_{n \to \infty} u_n - \lim_{n \to \infty} m_n = \lim_{n \to \infty} 2^{-n}(u_0-m_0)=0
\end{equation*}
So, $\displaystyle{\lim_{n \to \infty} u_n = \lim_{n \to \infty} m_n}$. If we put
\begin{equation*}
    M = \lim_{n \to \infty} u_n = \lim_{n \to \infty} m_n
\end{equation*}
then, by The Squeeze Theorem, we observe that $M$ is the maximum value of the function $f$ on the interval $[a,b]$.

Suppose that, at a certain stage in the process, the interval $[m_n,u_n]$ has just been defined. If the process is now stopped, the maximum value of $f$ is certain to lie in this interval. The best estimate of the maximum value at this stage is not $m_n$ or $u_n$ but the midpoint of the interval:
\begin{equation*}
    c_n=\frac{m_n+u_n}{2}
\end{equation*}
The error is then bounded as follows:
\begin{equation*}
    |M-c_n|\leq \frac{u_n-m_n}{2} =2^{-(n+1)}(u_0-m_0)
\end{equation*}
Summarizing this discussion, we have the following proposition.

{\bf Proposition 3.1.} If $[m_0,u_0 ],[m_1,u_1 ],\dots,[m_n,u_n ],\dots$ denote the intervals in the algorithm, then the limits $\displaystyle{\lim_{n \to \infty} u_n}$ and $\displaystyle{\lim_{n \to \infty} m_n}$ exist, are equal, and represent the maximum value of the function $f$ on the interval $[a,b]$. If $M=\displaystyle{\lim_{n \to \infty} c_n}$ and $c_n=\displaystyle{\frac{m_n+u_n}{2}}$, then
\begin{equation}\label{eq:2}
    |M-c_n|=2^{-(n+1)}(u_0-m_0)
\end{equation}
\begin{proof}
The proof is explained in the discussion above.
\end{proof}

{\bf Proposition 3.2.} The number of iterations needed, $n$, to achieve a given error (or tolerance), $\epsilon$, is given by 
\begin{equation*}
    n\geq \log_{2}{\frac{u_0-m_0}{\epsilon}}-1 = \frac{\log {(u_0-m_0)} - \log \epsilon}{\log 2} - 1
\end{equation*}
\begin{proof}
By Inequality [\ref{eq:2}], we can write
\begin{equation*}
    |M-c_n|=2^{-(n+1)}(u_0-m_0)
\end{equation*}
We want that $\displaystyle{|M-c_n|=2^{-(n+1)}(u_0-m_0)\leq \epsilon}$. If we solve the inequality $\displaystyle{2^{-(n+1)}(u_0-m_0 )\leq \epsilon}$ for $n$, we obtain that
\begin{equation*}
    n\geq \log_{2}{\frac{u_0-m_0}{\epsilon}}-1 = \frac{\log {(u_0-m_0)} - \log \epsilon}{\log 2} - 1
\end{equation*}
\end{proof}

{\bf Proposition 3.3.} The order of convergence of our algorithm is $1$. That is, the algorithm converges linearly.
\begin{proof}
The error at the $(n+1)$-th and $n$-th step is $|M-c_{n+1}|$ and $|M-c_n|$, respectively.
\begin{equation*}
    \lim_{n \to \infty} \frac{|M-c_{n+1}|}{|M-c_n|}=\frac{2^{-(n+2)}(u_0-m_0)}{{[2^{-(n+1)}(u_0-m_0)]}^p}=\frac{1}{2}{\left( \frac{2^n}{u_0-m_0} \right)}^{p-1}
\end{equation*}
Thus, we can conclude that the order of convergence $p$ is $1$ and error constant $C$ is $1$.
\end{proof}

{\bf Corollary 3.4.} The algorithm can be used to approximate the minimum value of a function $f$ which is continuous on the interval $[a,b]$.
\begin{proof}
If we approximate the maximum value of the function $y=-f(x)$ on the interval $[a,b]$ and change the sign of the approximated maximum value, then we get the minimum value of $f$.
\end{proof}

\section{Pros and Cons of the Algorithm}

There are a lot of algorithms and methods for unconstrained optimization. However, most of those methods and algorithms require differentiability criterion. Our algorithm does not require differentiability criterion. Thus, it can be applied any continuous function. What is more, we can remove continuity requirement because all we need is the Darboux Property which can be stated as
\begin{center}
    whenever $f(x_1)$ and $f(x_2)$ are different and $u$ is any number between them, then $f(x) = u$ for at least one $x$ between $x_1$ and $x_2$
\end{center}

Since continuous functions satisfy the Darboux property, we develop our algorithm over continuous functions. However, there are functions that are discontinuous at some points but the algorithm is still applicable. For example, consider the following piecewise function.

\[
  f(x) =
  \begin{cases}
                                   x+5 & \text{if \ $-4\leq x \leq -1$} \\
                                   4 & \text{if \ $-1 < x < 0$} \\
  3 & \text{if \ \ \ $0 \leq x \leq 1$}
  \end{cases}
\]

\begin{figure}[h]
    \centering
    \includegraphics[width=82mm,scale=0.5]{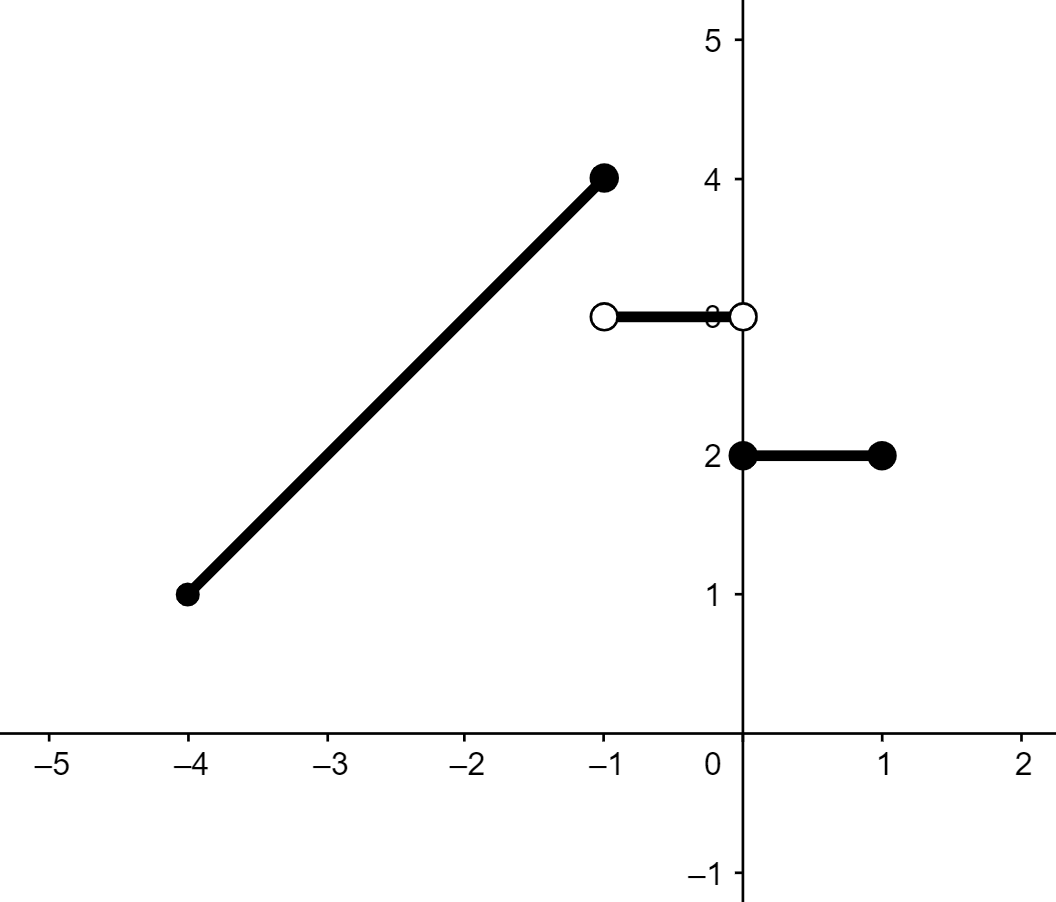}
    \caption{Graph of the piecewise function $f$}
    \label{graph}
\end{figure}

The function is clearly discontinuous at some points but it has the Darboux property. Therefore, our algorithm is still valid for the the piecewise function $f$.

Another advantage of the algorithm is that it always converges without a restriction on the initial starting points.

Now, let us have a look at some problematic aspects of the algorithm. The most restricting condition is that deciding whether the equation $f(x)=c$ has a solution or not may not be so easy. If we encounter a situation in which we cannot decide whether the equation $f(x)=c$ has a solution or not, then we cannot move on to the next step. 

Another issue is that there is no stopping criterion when we find the exact maximum value unless the given function is differentiable. Thus, we need to set an error to stop the algorithm. However, if the given function is differentiable, then we can decide that the algorithm find the exact maximum value. This topic will be discussed in the next section.

The last handicap is that the order of the convergence of the algorithm is $1$. Therefore, the algorithm converges, slowly.

\section{Analyzing Differentiable Functions}

Let $f$ be a differentiable function on the interval $[a,b]$. In the previous section, we emphasized that there is no stopping criterion when we find the exact maximum value, in general. However, if the function is differentiable, we can set a stopping criterion. If a function $f$ has a global or local maximum value at a point $x=c$, then $f'(c)$ must be equal to $0$. If $f(x)=c$ has a solution in the interval $[a,b]$ and we check whether $y=c$ is tangent to the function $f$, then we can decide whether we find the exact maximum value or not. 

\section{Further Works}
There are some further works to improve the algorithm. Since, I do not have enough time, I could not work on those subjects. However, anyone who is interested in this topic can study the following:

\begin{enumerate}
  \item It seems that the algorithm can be generalized into higher dimensions. Assume that $f:D\rightarrow \mathbb{R}$ is a function from $D\subseteq \mathbb{R}	^n$ to $\mathbb{R}$, where $D$ is a bounded, closed, and path connected subset of $\mathbb{R}	^n$. Then we can apply the Intermediate-Value Theorem and the Extreme Value Theorem for the function $f$. Thus, we may obtain an analogous version of the algorithm.
  \item The algorithm works for unconstrained optimization. However, it may be modified in such a way that the algorithm also works for constrained optimization.
  \item If the given function is a special class of functions(e.g. differentiable, convex, etc.), then the algorithm may be modified to work more efficiently.
  \item If someone needs just an integer value output, then the algorithm may be modified to produce the largest possible integer value output.
\end{enumerate}
\newpage
\section{References}

[1]. Kenneth A. Ross. \emph{Elementary Analysis}. 2nd Edition. 2013. Springer. New York. p 57.

[2]. George B. Thomas, Joel Hass, Christopher Heil, and Maurice D. Weir. \emph{Thomas' Calculus: Early Transcendentals}. 14th Edition. 2018. Pearson, Boston, p 70.

[3]. Serge Lang. \emph{Undergraduate Analysis}. 2nd Edition. 2005. Springer. New York. p 62.

[4]. Serge Lang. \emph{Undergraduate Analysis}. 2nd Edition. 2005. Springer. New York. p 61.

[5]. David Kincaid, Ward Cheney. \emph{Numerical Analysis: Mathematics of Scientific Computing}. 3rd Revised Edition. 2010. American Mathematical Society. Rhode Island. pp 74-78.

\newpage



\section{Appendix}
{\bf An example Matlab code which implements our algorithm:}
\begin{lstlisting}[frame=single]
syms x  
eqn = input(sprintf('Enter a continuous function as a function of x'));
a = input(sprintf('Enter the lower bound of the interval [a,b]'));
b = input(sprintf('Enter the upper bound of the interval [a,b]'));
if b<=a
    fprintf('Error: b must be bigger than a\n')
    return
end
u = input(sprintf('Enter a number u such that f(x)=u has a solution in the interval [a,b]\n'));
k = solve(eqn==u,x,'Real',true);
t = isempty(k);
g = max(size(k));
z = 0;
if t == 1 
    fprintf('Error: Enter a number u such that f(x)=u has a solutionin the interval [a,b]\n') 
  return
end
for i=1;g;
    if k(i)>a && k(i)<b
       z=1;
    end
end
if z==0
    fprintf('Error: Enter a number u such that f(x)=u has a solutionin the interval [a,b]\n')
    return
end
v = input(sprintf('Enter a number v bigger than u such that f(x)=v has no solution'));
l = solve(eqn==v,x,'Real',true);
r = isempty(l);
o = 0;
for i=1;g;
    if r==0
        if l(i)>a && l(i)<b
            o=1;
        end
    end
end
if r == 0 && o==1
    fprintf('Error: Enter a number v bigger than u such that f(x)=v has no solution')
    return
end
if v<=u
    fprintf('Error: Enter a number u such that f(x)=u has a solutionin the interval [a,b]\n')
    return
end
nmax = input(sprintf('Enter the maximum number of iterations'));
e = input(sprintf('Enter the maximum error'));
if nmax<0 || e<0
    fprintf('Error: Enter a positive number for the number of iterations and the error')
    return
end
n=1;
while n<=nmax
    p=(u+v)/2;
    k = solve(eqn==p,x,'Real',true);
    t = isempty(k);
    g = max(size(k));
    z=0;
    if t == 1 
    v=p;
    end
    for i=1;g;
        if t==0
            if k(i)>a && k(i)<b       
                z=1;
            end
         end
    end
    if z==1
       u=p;
    else 
       v=p;
    end
    if v-u<e
        fprintf('approximate maximum is %f',(u+v)/2)
        return
    end
    n=n+1;
    if n>nmax
        fprintf('Error: please enter a bigger number of maximum iteration to approximate the maximum with an error less than %f',e)
    end
end
\end{lstlisting}

\end{document}